\documentclass[11pt]{article}
\usepackage{graphicx}
\usepackage{amsmath}
\usepackage[dvips]{epsfig}
\usepackage{amssymb}
\begin{document}
\begin{center} {\Large \bf  On a $p$-adic interpolation function for the multiple generalized Euler numbers and its derivatives }
\\ \vspace*{12 true pt}  Taekyun   Kim
\vspace*{12 true pt} \\
 Division of General Education-Mathematics,\\
 Kwangwoon University, Seoul 139-701,  Korea \\
 e-mail: tkkim@kw.ac.kr
                        \end{center}
%--------------------------------------------------------------------
\vspace*{12 true pt} \noindent {\bf Abstract :} We study analytic
function interpolating the multiple generalized Euler numbers
attached to $\chi$ at  negative integers in complex plane and  we
consider the multiple $p$-adic $l$-function as the $p$-adic analog
of the above function. Finally, we give the value of the partial
derivative of this multiple $p$-adic $l$-function at $s=0$.

\vspace*{12 true pt} \noindent {\bf 2000 Mathematics Subject
Classification :} 11B68, 11S40, 11S80

 \vspace*{12 true pt} \noindent {\bf
Key words :}   Euler numbers, Euler polynomials,  multiple
generalized Euler numbers, multiple generalized Euler numbers
attached to $\chi$, multiple $p$-adic $l$-function

\begin{center} {\bf 1. Introduction } \end{center}
 Let $p$ be a fixed odd prime number. Throughout this paper  $\mathbb{Z}_p, \mathbb{Q}_p, \mathbb{C}$  and $\mathbb{C}_p$
 will, respectively, denote the ring of $p$-adic rational integers,
  the field of $p$-adic rational numbers,   the complex number
  field and the  completion of algebraic closure of   $\mathbb{Q}_p$.
Let $\mathbb{N}$ be   the  set of natural numbers and let $\nu_p$
be the normalized exponential valuation of $\mathbb{C}_p$ with
$|p|_p=p^{-\nu_p(p)}=p^{-1}.$ The Euler numbers in $\mathbb{C}$
are defined by the formula
$$F(t)=\dfrac{2}{e^t+1}=\sum_{n=0}^\infty E_n \dfrac{t^n}{n!} \text{ for } |t|< \pi, \text{ see [1-6]}. \eqno(1) $$
It follows from the definition that $ E_0=1, E_1=-1/2, E_2=0,
E_3=1/4, \cdots, $ and $E_{2k}=0$ for $k=1,2, 3, \cdots.$

Let $ r \in \mathbb{N}$. Then the multiple Euler numbers of order
$r$ are defined as
$$F^{(r)}(t)= \left( \dfrac{2}{e^t+1} \right)^{r} =\sum_{n=0}^\infty E_n^{(r)} \dfrac{t^n}{n!},  |t|< \pi . $$
Let $x$ be a variable. Then the multiple Euler polynomials are
also defined by the rule
$$F^{(r)}(t,x)= \left( \dfrac{2}{e^t+1} \right)^{r}e^{xt} =\sum_{n=0}^\infty E_n^{(r)}(x) \dfrac{t^n}{n!} . \eqno(2) $$

For  $f \in \mathbb{N}$ with $ f \equiv 1 \pmod{2}$, assume that
$\chi$ is a  primitive  Dirichlet character with conductor $f$. It
is known that  the generalized Euler numbers  attached to $\chi$,
$E_n, \chi$, are defined by the rule
$$  F_{\chi} (t)= 2 \sum_{a=1}^{ f} \dfrac{\chi(a) (-1)^a e^{at}}{e^{ft} +1 }  = \sum_{n=0}^\infty  E_{n, \chi}
   \dfrac{t^n}{n!}, \eqno(3)  $$
   where $|t|< \frac{\pi}{f},$ (see [1]).

In this paper, we consider the multiple generalized Euler numbers
attached to $\chi$ in the sense of the multiple of Eq.(3). From
these numbers, we study an analytic function interpolating the
multiple generalized Euler numbers attached to $\chi$ at  negative
integers in complex number field. In the sense of $p$-adic analog
of the above function and give the values of the partial
derivative of this multiple  $p$-adic $l$-function at $s=0$.

\bigskip

\begin{center} {\bf 2. Analytic functions associated with Euler numbers and polynomials}
\end{center}

We now consider the multiple generalized Euler numbers attached to
$\chi$, $E_{n, \chi}^{(r)}$, which are defined by
$$F_{\chi}^{(r)}(t)=2^r \sum_{a_1, \cdots, a_r=1}^f \dfrac{(-1)^{\sum_{i=1}^r a_i} \chi{(a_1+ \cdots + a_r)} e^{t\sum_{i=1}^r a_i}}{(e^{ft}+1)^r}
=\sum_{n=0}^\infty E_{n, \chi}^{(r)} \dfrac{t^n}{n!}. \eqno(4)$$
By (2) and (4), we readily see that
$$E_{n, \chi}^{(r)}
=f^n \sum_{a_1, \cdots, a_r=0}^{f-1}  (-1)^{\sum_{i=1}^r a_i}
\chi{(a_1+ \cdots + a_r)} E_n^{(r)} \left( \dfrac{a_1+ \cdots +
a_r}{f} \right) . \eqno(5)$$
 For $ s \in \mathbb{C}$, we have
$$2^r \sum_{n_1, \cdots, n_r=0}^\infty  \dfrac{(-1)^{n_1+ \cdots +
n_r}}{(x+n_1+ \cdots +n_r)^s}= \dfrac{1}{\Gamma(s)}\int_0^\infty
F^{(r)}(-t, x) t^{s-1} dt, \eqno(6)$$ where $x \neq 0, -1, -2,
\cdots.$ From (6), we can consider the multiple Euler zeta
function as follows:
$$\zeta_r(s, x)=2^r \sum_{n_1, \cdots, n_r=0}^\infty  \dfrac{(-1)^{n_1+ \cdots +
n_r}}{(x+n_1+ \cdots +n_r)^s}, \text{ for } x\in \mathbb{C}, x
\neq 0, -1, -2, \cdots. \eqno(7)$$ By (1), (2), (5), (6) and (7),
we easily see that
$$\zeta_r(-n, x)=E_n^{(r)}(x) \text{ for } n \in \mathbb{N}.$$
By using complex integral and (4), we can also obtain the
following equation: for $s \in \mathbb{C}$,
$$ \dfrac{1}{\Gamma(s)}\int_0^\infty
F_{\chi}^{(r)}(-t) t^{s-1} dt= 2^r \sum_{\begin{subarray}{l} n_1,
\cdots, n_r=0 \\ n_1+ \cdots+ n_r \neq 0 \end{subarray}}^\infty
\dfrac{{\chi (n_1+ \cdots +n_r)} (-1)^{n_1+ \cdots + n_r }}{ (n_1+
\cdots +n_r)^s}, \eqno(8)$$ where  $\chi$ is the primitive
Dirichlet  character with conductor $f \in \mathbb{N}$ with  $ f
\equiv 1 \pmod{2}$.

By (8), we define Dirichlet's type multiple Euler $l$-function in
complex plane as follows: for $s \in \mathbb{C}$,
$$ l_r(s, \chi)= 2^r \sum_{\begin{subarray}{l} n_1,
\cdots, n_r=0 \\ n_1+ \cdots+ n_r \neq 0 \end{subarray}}^\infty
\dfrac{{\chi (n_1+ \cdots +n_r)} (-1)^{n_1+ \cdots + n_r }}{ (n_1+
\cdots +n_r)^s}. \eqno(9)$$ From (4), (8) and (9), we note that
$$l_r(-n, \chi) = E_{n, \chi}^{(r)} \text{ for } n \in \mathbb{N}. \eqno(10)$$
Let $s$ be a complex variable, and let $a$ and $F$ be integer with
$0< a< F$ and $ F \equiv 1 \pmod{2}$. Then we consider partial
zeta function $T_r(s; a_1, \cdots, a_r | F)$ as follows:
$$ \aligned T_r(s; a_1, \cdots, a_r | F) & =2^r \sum_{\begin{subarray}{l} m_1,
\cdots, m_r>0 \\ m_i \equiv a_i\pmod{F}
\end{subarray}} \dfrac{ (-1)^{m_1+
\cdots + m_r }}{ (m_1+ \cdots +m_r)^s} \\
&=(-1)^{a_1+\cdots+a_r} F^{-s} \zeta_r \left( s,
\dfrac{a_1+\cdots+a_r}{F} \right). \endaligned \eqno(11)
 $$
Let  $\chi( \neq 1)$ be  the  Dirichlet  character with conductor
$f \in \mathbb{N}$ with  $ f \equiv 1 \pmod{2}$ and let $F$ be the
multiple of $f$ with $ F \equiv 1 \pmod{2}$. Then Dirichlet's type
multiple  $l$-function can be expressed  as the sum
$$ l_r(s, \chi)= \sum_{a_1, \cdots, a_r=1}^F \chi(a_1+\cdots+a_r)  T_r(s; a_1, \cdots, a_r | F)
\text{ for } s \in \mathbb{C}. \eqno(12)$$ By simple calculation,
we easily see that
$$T_r(-n; a_1, \cdots, a_r | F)= F^n (-1)^{a_1+ \cdots + a_r } E_n^{(r)} \left(
\dfrac{a_1+\cdots+a_r}{F} \right) \text{ for } r, n \in
\mathbb{N}. \eqno(13)$$ From (5), (12) and (13), we note that
$$l_r(-n, \chi)=E_{n, \chi}^{(r)}  \text{ for }  n \in
\mathbb{N}. \eqno(14)$$ By (1) and (2), we have
$$E_n^{(r)}(x)= \sum_{l=0}^n \binom nl E_l^{(r)} x^{n-l} = \sum_{l=0}^n \binom nl E_{n-l}^{(r)} x^{l}. \eqno(15)$$
From (11), (13) and (15), we derive the following equation:

$$T_r(s; a_1, \cdots, a_r | F)= (-1)^{a_1+ \cdots + a_r }(a_1+ \cdots + a_r)^{-s} \sum_{k=0}^\infty \binom{-s}{k}  \left(
\dfrac{F}{a_1+\cdots+a_r} \right)^k E_k^{(r)}. \eqno(16)$$ It
follows from (12), (13) and (16) that
$$ \aligned  l_r(s, \chi)   = &
\sum_{a_1, \cdots, a_r=1}^F \chi(a_1+\cdots+a_r)  (-1)^{a_1+ \cdots + a_r }(a_1+ \cdots + a_r)^{-s} \\
& \times \sum_{m=0}^\infty \binom{-s}{m} \left(
\dfrac{F}{a_1+\cdots+a_r} \right)^m E_m^{(r)}. \endaligned
\eqno(17)
 $$
 The values of $l_r(s, \chi)$ at negative integer are algebraic
 numbers, and hence they can be regarded as numbers belonging to
 an extension of $ \mathbb{Q}_p$.  In the next section, we
 therefore look for a $p$-adic function that agrees with $l_r(s,
 \chi)$ at negative integers.

\bigskip

\begin{center} {\bf 3. On a $p$-adic interpolation function for the multiple Euler numbers and its derivative}
\end{center}

In this section, we study the $p$-adic analogs of the Dirichlet's
type multiple  Euler $l$-function, $ l_r(s, \chi)$, which were
introduced in previous section. If fact, this function is the
$p$-adic interpolation function for the multiple generalized Euler
numbers attached to $\chi$ at negative integers. Let $w$ denote
the Teichm\"{u}ller character with the conductor $f_w=p$. For an
arbitrary character $\chi$, we set $\chi_n = \chi w^{-n}, n \in
\mathbb{Z},$  in the sense of the product of characters. Let
$$<a>=w^{-1}(a)a =\dfrac{a}{w(a)}.$$
Then we note that $ <a> \equiv 1 \pmod{ p \mathbb{Z}_p}$. Let
$$ A_j(x)= \sum_{n=0}^\infty a_{n,j} x^n,  a_{n,j} \in
\mathbb{C}_p, j=0,1,2, \cdots $$ be a sequence of power series,
each convergent on a fixed subset
$$ D= \{ s \in \mathbb{C}_p \mid |s|_p < p^{1-\frac{1}{p-1} } \},$$
of $\mathbb{C}_p$ such that
$$ \aligned  (1) & \quad    a_{n,j} \rightarrow a_{n,0} \text{ as } j \rightarrow \infty  \text{ for any } n ; \\
 (2) &\quad   \text{ for each } s \in D \text{ and } \epsilon
>0, \text{ there exists an  } n_0= n_0(s, \epsilon) \text{ such that }  \\
& \quad  \quad \quad \quad \quad  \quad \quad \quad | \sum_{n \geq
n_0} a_{n,j} s^n |_p < \epsilon \text{  for  } \forall j.
\endaligned
 $$
In this case,
$$\lim_{j \rightarrow \infty} A_j(s)=A_0(s), \text{ for all } s \in D.$$
This was used by Washington ([6]) to show that each of functions
$w^{-s}(a) a^s$ and $$ \sum_{m=0}^\infty \binom sm \left(
\dfrac{F}{a}\right)^m B_m, $$ where $F$ is multiple of $p$ and $f$
and $B_m$  is the $m$-th Bernoulli numbers, is analytic on $D$
(see [6]).

 Let  $\chi$ be  a primitive Dirichlet's
character with conductor $f \in \mathbb{N}$ with  $ f \equiv 1
\pmod{2}$. Then we consider the multiple Euler $p$-adic
$l$-function, $l_{p,r}(s, \chi)$, which interpolates the multiple
generalized Euler numbers attached to $\chi$ at negative integers.

For $f \in \mathbb{N}$ with  $ f \equiv 1 \pmod{2}$, let us assume
that $F$ is a positive integral multiple of $p$ and $f=f_\chi$. We
now define the multiple Euler $p$-adic $l$-function as follows:
$$ \aligned  l_{p,r}(s, \chi)   = &
\sum_{a_1, \cdots, a_r=1}^F \chi(a_1+\cdots+a_r)<a_1+ \cdots + a_r>^{-s} (-1)^{a_1+ \cdots + a_r } \\
& \times  \sum_{m=0}^\infty \binom{-s}{m} \left(
\dfrac{F}{a_1+\cdots+a_r} \right)^m E_m^{(r)}. \endaligned
\eqno(18)
 $$
From (18), we note that $ l_{p,r}(s, \chi) $ is analytic for $s
\in D$.

For $n \in \mathbb{N},$ we have
$$E_{n, \chi_n}^{(r)}
=F^n \sum_{a_1, \cdots, a_r=0}^{F-1}  E_n^{(r)} \left( \dfrac{a_1+
\cdots + a_r}{F} \right) \chi_n{(a_1+ \cdots + a_r)} (-1)^{a_1+
\cdots + a_r} . \eqno(19)$$ If $\chi_n(p) \neq 0$, then $(p,
f_{\chi_n})=1,$ and thus the ratio $F/p$ is a multiple of
$f_{\chi_n}$.

Let $$ I_0 = \left \{ \dfrac{a_1+ \cdots + a_r}{p} \mid  a_1+
\cdots + a_r \equiv 0\pmod{ p} \text{ for some } a_i \in
\mathbb{Z} \text{ with } 0 \leq a_i \leq F \right \}.$$ Then we
have
$$ \aligned & F^n \sum_{\begin{subarray}{l} a_1,
\cdots, a_r=0 \\ p \mid a_1+ \cdots +a_r
\end{subarray}}^{F-1} E_n^{(r)} \left(
\dfrac{a_1+\cdots+a_r}{F} \right) \chi_n (a_1+\cdots+a_r)
(-1)^{a_1+ \cdots + a_r } \\
& = p^n \left( \dfrac{F}{p} \right)^n \chi_n(p)
\sum_{\begin{subarray}{l} a_1, \cdots, a_r=0 \\ \quad \beta \in
I_0
\end{subarray}}^{F/p} \chi_n (\beta) (-1)^\beta  E_n^{(r)} \left( \dfrac{\beta}{F/p}
\right).
\endaligned
\eqno(20)
 $$
From (20), we define the second multiple generalized Euler numbers
attached to $\chi$ as follows:
$$E_{n, \chi_n}^{* (r)}= \left( \dfrac{F}{p} \right)^n \sum_{\begin{subarray}{l} a_1, \cdots, a_r=0 \\ \quad \beta \in I_0
\end{subarray}}^{F/p} \chi_n (\beta) (-1)^\beta  E_n^{(r)} \left( \dfrac{\beta}{F/p}\right). \eqno(21)$$

By (19), (20) and (21), we easily see that
$$ \aligned &  E_{n, \chi_n}^{ (r)}- p^n \chi_n(p) E_{n, \chi_n}^{* (r)} \\
&= F^n \sum_{\begin{subarray}{l} a_1, \cdots, a_r=1 \\ p \nmid
a_1+ \cdots +a_r\end{subarray}}^F   \chi_n (a_1+\cdots+a_r)
(-1)^{a_1+ \cdots + a_r }  E_n^{(r)} \left(
\dfrac{a_1+\cdots+a_r}{F} \right).
\endaligned
\eqno(22)
 $$
 By the definition of the multiple Euler polynomials of order $r$,
 we see that
$$ E_n^{(r)} \left(
\dfrac{a_1+\cdots+a_r}{F} \right)=F^{-n}(a_1+\cdots+a_r)^n
\sum_{k=0}^n \binom{n}{k} \left( \dfrac{F}{a_1+\cdots+a_r}
\right)^k E_k^{(r)}.\eqno(23)$$ From (22) and (23), we have

$$ \aligned  E_{n, \chi_n}^{ (r)}- p^n \chi_n(p) E_{n, \chi_n}^{* (r)} & = \sum_{\begin{subarray}{l} a_1, \cdots, a_r=1 \\ p \nmid a_1+
\cdots +a_r\end{subarray}}^F  (a_1+\cdots+a_r)^n  \chi_n
(a_1+\cdots+a_r)
 (-1)^{a_1+ \cdots + a_r } \\
& \quad  \quad \quad  \quad  \times \sum_{k=0}^n \binom{n}{k}
\left( \dfrac{F}{a_1+\cdots+a_r} \right)^k E_k^{(r)} .
\endaligned
\eqno(24)
 $$
By (18) and (24), we readily see that
$$ \aligned & l_{p, r}(-n, \chi) \\
 & =
\sum_{\begin{subarray}{l} a_1, \cdots, a_r=1 \\ p \nmid a_1+
\cdots +a_r\end{subarray}}^F  (a_1+\cdots+a_r)^n  \chi_n
(a_1+\cdots+a_r)
 (-1)^{a_1+ \cdots + a_r } \\
& \quad  \quad  \quad  \quad  \times \sum_{m=0}^n \binom{n}{m}
\left(\dfrac{F}{a_1+\cdots+a_r} \right)^m E_m^{(r)} \\
&=E_{n, \chi_n}^{ (r)}- p^n \chi_n(p) E_{n, \chi_n}^{* (r)}.
\endaligned
\eqno(25)
 $$
Therefore, we obtain the following theorem.

\bigskip
{ \bf Theorem 1.} Let $F$ be a positive integral of $p$ and
$f(=f_{\chi_n})$, and let
$$ \aligned  l_{p,r}(s, \chi)   = &
\sum_{a_1, \cdots, a_r=1}^F \chi(a_1+\cdots+a_r)<a_1+ \cdots + a_r>^{-s} (-1)^{a_1+ \cdots + a_r } \\
& \times  \sum_{m=0}^\infty \binom{-s}{m} \left(
\dfrac{F}{a_1+\cdots+a_r} \right)^m E_m^{(r)}. \endaligned
 $$
 Then $ l_{p,r}(s, \chi)$ is analytic on $D$.
 Furthermore, for each $n \in \mathbb{N}$, we have
 $$ l_{p,r}(-n, \chi)=E_{n, \chi_n}^{ (r)}- p^n \chi_n(p) E_{n, \chi_n}^{* (r)}.$$
\bigskip

Using Taylor expansion at $s=0$, we get
$$\binom{-s}{m}=\dfrac{(-1)^m}{m}s+ \cdots  \text{ if } m \geq 1. \eqno(26)$$
From (26) and Theorem 1, we obtain the following corollary.

\bigskip
{ \bf Corollary 2.} Let $F$ be a positive integral multiple  of
$p$ and $f$. Then we have
$$ \aligned  \dfrac{\partial }{\partial s}l_{p,r}(0, \chi)   = &
 \sum_{\begin{subarray}{l} \quad a_1, \cdots, a_r=1 \\ (a_1+
\cdots+a_r, p)=1\end{subarray}}^F \chi(a_1+\cdots+a_r)(-1)^{a_1+ \cdots + a_r }\left(1- \log_p(a_1+ \cdots+a_r)\right) \\
& +  \sum_{\begin{subarray}{l} \quad a_1, \cdots, a_r=1 \\ (a_1+
\cdots+a_r, p)=1\end{subarray}}^F \chi(a_1+\cdots+a_r)(-1)^{a_1+
\cdots + a_r } \sum_{m=1}^\infty \dfrac{(-1)^m}{m} \left(
\dfrac{F}{a_1+\cdots+a_r} \right)^m E_m^{(r)}, \endaligned
 $$ where $\log_p x$ is denoted by the $p$-adic logarithm.

\bigskip

\bigskip
\begin{center}{\bf REFERENCES}\end{center}
\begin{enumerate}

\item
{T. Kim,} { Barnes type multiple $q$-zeta function and $q$-Euler
polynomials,} { J. Phys. A: Math. Theor.}, {43}(2010)
255201(11pp).

\item
{ T. Kim,}  {   $q$-Volkenborn integration}, { Russ. J. Math.
phys.},  {  9}(2002), 288-299.

\item
{T. Kim,}  {  Note on the Euler $q$-zeta functions,} {J.  Number
Theory}, { 129}(2009), 1798-1804.

\item
{H. Ozden, Y. Simsek, S.-H. Rim, I. N. Cancul,}  { A note on
$p$-adic $q$-Euler measure,} { Adv. Stud. Contemp. Math.},
{14}(2007), 233-239.

\item
{ K. Shiratani, S. Yamamoto}  {   On a $p$-adic interpolation
function for the Euler numbers and its derivatives}, { Mem. Fac.
Sci. Kyushu Univ. A}, {39}(1985), 113-125.

\item
{L. C. Washington,}  {  Introduction to Cyclotomic Field,}
{Springer}, 1982.

\end{enumerate}

\end{document}